\documentclass[11pt] {article} 
  \usepackage{amsmath}
    \usepackage{amssymb}

\newtheorem{theorem}{Theorem}[section]
\newtheorem{lemma}{Lemma}[section]

\newtheorem{remark}{Remark}[section]

\newcommand{\eqnsection}{
   \renewcommand{\theequation}{\thesection.\arabic{equation}}
   \makeatletter
   \csname @addtoreset\endcsname{equation}{section} 
   \makeatother}

\def \ov{\overline}

\def \be{\begin{equation}}
\def \ee{\end{equation}}
\def \bt{\begin{theorem}} 
\def \et{\end{theorem}}
\def \bl{\begin{lemma}} 
\def \el{\end{lemma}}
\def \bea{\begin{eqnarray}}
\def \eea{\end{eqnarray}}
\def \bas{\begin{eqnarray*}}
\def \eas{\end{eqnarray*}}

\def \de{\delta}
\def \De{\Delta}

\def \vf{\varphi}
\def \si{\sigma}

\def \ff{\infty}

\def \wt{\widetilde}

\def\stl{\stackrel{law}{=}}

\def \AA{{\cal A}}

\def \HH{{\cal H}}

\def\b1{\mathbf 1}
\def \({\left(}
\def \){\right)}

\def \nn{\nonumber}

\def \bc{\begin{center} }
\def \ec{\end{center} }
\def \bs{\begin{slide} }
\def \es{\end{slide} }

\def\square{{\vcenter{\vbox{\hrule height.3pt
        \hbox{\vrule width.3pt height5pt \kern5pt
           \vrule width.3pt}
        \hrule height.3pt}}}}
\def\qed{{\hfill $\square$ \bigskip}}

\eqnsection

\begin{document}

\title{Local and  uniform moduli of continuity of chi--square processes }

  \author{  Michael B. Marcus\,\, \,\, Jay Rosen \thanks{Research of     Jay Rosen was partially supported by  grants from the Simons Foundation.   }}
\maketitle
 \footnotetext{ Key words and phrases:  chi--square processes, moduli of continuity }
 \footnotetext{  AMS 2020 subject classification:  60F15,  60G15, 60G17}

\maketitle
 
\section{Introduction}

  Let $\eta=\{\eta(t);t\in [0,1]\}$ be a mean zero    continuous Gaussian process with covariance $U=\{U(s,t),s,t\in [ 0,1]\},$ with $U(0,0)>0$.  Let $\{\eta_{i};i=1,\ldots, k\}$ be independent copies of $\eta$  and set,
  \begin{equation} \label{1.9mm}
  Y_{k}(t)=\sum_{i=1}^{k} \eta^2_{i}(t),\qquad t\in [ 0,1].\end{equation}
The stochastic process  $Y_{k } =\{Y_{k }(t),t\in  [ 0,1] \}$ is referred to as a chi--square process of order $k $  with kernel $U$.

   \bt \label{theo-1.3m}   Let $\phi(t)$   be a positive function on $[0,\de]$ for  some $\de>0$.  
 If
 \be\limsup_{t\to 0}\frac{  \eta(t)-\eta(0)}{ \phi(t)   }=1 \qquad 
a.s.\label{rev.1wq},
\ee
then for all integers $k\ge 1$,
\begin{equation} 
 \limsup_{t\to 0} \frac{Y_{k }(t)-Y_{k }(0)} {   \phi  (t)} =  2 Y^{1/2}_{k}(0) \qquad a.s.\label{rev.2qq}
\end{equation} 
 
 \et
 
   When $k=1$ this is particularly simple. Since $\eta$ is symmetric it follows from (\ref{rev.1wq}) that,
 \be
\liminf_{t\to 0}\frac{  \eta(t)-\eta(0)}{ \phi(t)   }=-1 \qquad 
a.s.\label{rev.1wqb}
\ee
Therefore, writing $ Y_{1 }(t)-Y_{1}(0)=(\eta  (t)-\eta  (0))(\eta  (t)+\eta  (0))$ and using the continuity of   $\eta$, we see that 
\begin{equation}
 \limsup_{t\to 0} \frac{Y_{1 }(t)-Y_{1 }(0)  }  {   \phi  (t)} =  2\(\eta(0)\vee -\eta(0)\)  \qquad a.s.,\label{rev.2ff}
\end{equation} 
which is (\ref{rev.2qq}). 

 A result similar to Theorem \ref{theo-1.3m}
for the limiting behavior of chi--square sequences at infinity is given in  \cite[Lemma 6.5]{asym}.

\medskip   Set  \be \si^2(u,v)=E(\eta(u)-\eta(v))^2\quad\text{and}\quad \wt\si^2(x)=\sup_{|u-v|\le x}\si^2(u,v).\label{1.15}
\end{equation}

\bt  \label{theo-2.1m} Assume that  $\inf_{t\in [0,1]}U(t,t)>0$ and, 
\begin{equation}
\lim_{x\to0}\wt\si^2(x)\log 1/x
=0.\label{ovabs.3}\end{equation}
Let  $\vf(t)$   be a positive function on $[0,1]$. Then if 
 \be
\lim_{h\to 0}\sup_{\stackrel{|u-v|\le h }{ u,v\in\De}}\frac{  \eta(u)-\eta(v)}{ \vf(|u-v|)   }=1 \qquad 
a.s.\label{rev.1wqu},
\ee
 for all intervals $\De\subset [0,1]$,  
it follows that for all intervals $\De\subset [0,1]$ and   all  integers  $k\ge 1$,  
\begin{equation} \lim_{h\to 0}\sup_{\stackrel{|u-v|\le h }{ u,v\in\De}}  \frac{Y_{k }(u)-Y_{k }(v) }{  \vf  (|u-v|)} =   2 \sup_{u\in\De}Y_{k }^{1/2}(u), \hspace{.2 in}a.s.\label{rev.2quu}
\end{equation} 
 
 \et

  When $\eta$ is a continuous Gaussian process with stationary increments, $ \si^2(u,v)$ in (\ref{1.15}) can be written as 
 $ \si^2(u-v,0)$. In this case if $ \wt \si^2(x)$ is asymptotic to an increasing function at 0, 
 then (\ref{rev.1wqu}) implies (\ref{ovabs.3}). We discuss this further in Remark \ref{rem-2.1}.
 
  \medskip   An extensive treatment of Gaussian processes satisfying (\ref{rev.1wq}) and (\ref{rev.1wqu})
is given  in \cite[Chapter 7]{book}.

\section{Proofs}
  
 \medskip
\noindent{\bf  Proof of Theorem \ref{theo-1.3m} } Let $\eta_{i}(t)$, $i=1,\ldots,k$, be independent copies of $\eta(t)$.   
  We write, 
\begin{eqnarray}
 \eta_{i}^{2}(t)- \eta_{i}^{2}(0)&=&( \eta_{i}(t)- \eta_{i}(0))( \eta_{i}(t)+ \eta_{i}(0))
\label{16.40}\\
&=& ( \eta_{i}(t)- \eta_{i}(0))( 2 \eta_{i}(0)+( \eta_{i}(t)- \eta_{i}(0)))   \nonumber\\
&=& 2( \eta_{i}(t)- \eta_{i}(0))  \eta_{i}(0)+( \eta_{i}(t)- \eta_{i}(0))^{2}.  \nonumber
\end{eqnarray}
 By (\ref{rev.1wq})
 \bea
 &&  \limsup_{t\to 0}\frac{\sum_{i=1}^{k}( \eta _{i }(t)-\eta_{i }(0))^{2} }{  \phi (t)}\label{5.8mm}\\
 &&\nn\qquad\le \sum_{i=1}^{k} \limsup_{t\to 0}\frac{| \eta _{i }(t)-\eta_{i }(0)| }{  \phi (t)}\lim_{t\to0}| \eta _{i }(t)-\eta_{i }(0)|=0.
   \eea
Consequently, using (\ref{1.9mm}) we see that,  
\begin{equation}
\limsup_{t\to0}\frac{Y_{k }(t)-Y_{k }(0)}{   \phi (t)}=\limsup_{t\to0}\frac{2\sum_{i=1}^{k}( \eta _{i }(t)-\eta _{i }(0))\eta _{i }(0) }{   \phi (t)}\label{5.4mmt}.
   \end{equation}
Write,  
\bea
  &&( \eta_{i}(t)- \eta_{i}(0))  \eta_{i}(0) \label{16.41}\\
  &&\qquad=\( \eta_{i}(t)-\frac{U(0,t)}{U(0,0)} \eta_{i}(0)\)  \eta_{i}(0) -\(\frac{U(0,0)-U(0,t)}{U(0,0)}\)\eta^{2}_{i}(0)\nn.\label{2.5mm}
   \eea
  We show below that  
  \begin{equation}
\limsup_{t\to0}\frac{|U(0,0)-U(0,t)|}{U(0,0)\phi (t)}=0.\label{2.30mm}
   \end{equation}
   Consequently,
   \begin{equation}
\limsup_{t\to0}\frac{Y_{k }(t)-Y_{k }(0)}{  \phi (t)}=\limsup_{t\to0}\frac{2\sum_{i=1}^{k}\( \eta_{i}(t)-\frac{U(0,t)}{U(0,0)} \eta_{i}(0)\)  \eta_{i}(0)}{  \phi (t)}\label{5.4mm}.
   \end{equation}

 	Let $\{\xi_{i}(t),t\in [0,1]\}$, $i=1,\ldots,k$, be independent copies of a mean zero Gaussian process  $\{\xi (t),t\in [0,1]\}$, and set $\vec \xi(t)=(\xi_{1}(t),\ldots,\xi_{k}(t))$. Let $\vec v\in R^{k}$ with $\|\vec v\|_{2}=1$. Computing the covariances we see that,
\begin{equation}
  \{ (\vec v\cdot \vec\xi(t)),t\in [0,1]\}\stl \{\xi(t),t\in [0,1]\}\label{2.32mm}.
   \end{equation}
  (This relationship is used by P. Revesz in  \cite[Theorem 18.1]{Rev} to obtain LILs for   Brownian motion in $R^{k}$.)

Therefore, since 
   $\( \eta_{i}(t)-(U(0,t)/U(0,0)) \eta_{i}(0)\)$  and $\eta_{i}(0)$ are independent for $i=1,\ldots,k$, we see that,
   \bea
&& \Big  \{   \(\vec  \eta (t)-\frac{U(0,t)}{U(0,0)} \vec\eta (0)\)\cdot  \frac{\vec\eta (0)}{ \|\vec\eta (0)\|_{2}},t\in [0,1]\Big\}\label{2.10mm}\\
&&\qquad\stl  \Big  \{   \(   \eta (t)-\frac{U(0,t)}{U(0,0)} \eta (0)\),t\in [0,1]\Big\}\nn,
      \eea
where  $\vec \eta(t)=(\eta_{1}(t),\ldots,\eta_{k}(t))$.
  Consequently,  (\ref{5.4mm}) implies that 
\begin{equation}
\limsup_{t\to0}\frac{Y_{k }(t)-Y_{k }(0)}{ \phi (t)  \|\vec\eta (0)\|_{2}}\stl   \limsup_{t\to0}\frac{ 2 \( \eta (t)-\frac{U(0,t)}{U(0,0)} \eta (0)\)  }{  \phi (t)}.
   \end{equation}
   Using (\ref{16.41}) again and (\ref{2.30mm}) we see that this  implies that, 
   \be
   \limsup_{t\to0}\frac{Y_{k }(t)-Y_{k }(0)}{    \phi (t) \|\vec\eta (0)\|_{2}}\stl    \limsup_{t\to0}\frac{ 2(  \eta (t)- \eta (0) )   }{   \phi (t)}=2,
   \ee
   where the last equality uses (\ref{rev.1wq}). Since,
\begin{equation}
  \|\vec \eta  (0)\|_{2}=\(\sum_{i=1}^{k} \eta^{2}_{i}  (0) \)^{1/2}=  Y_{k }^{1/2}(0),\label{2.20mm}
\end{equation}
  we obtain (\ref{rev.2qq}).

  \medskip  To obtain    (\ref{2.30mm}) we first note that  it follows from (\ref{rev.1wq}) that, 
   \begin{equation}
   \phi(t)= (E(\eta (t)-  \eta (0))^2)^{1/2}h(t),\label{c2.12}
   \end{equation}
   for some function $h$ such that $\lim_{t\downarrow0}h(t)=\ff $.   To see this, suppose that it is false. Then there exists a sequence $\{t_k \}$, with $\lim_{k\to\ff}t_k=0$, such that $\sup_k h(t_k)\le M$. Therefore, if (\ref{rev.1wq}) holds, we would have,
 \be
\sup_{k}\frac{  \eta(t_k)-\eta(0)}{ (E(\eta (t_{k})-  \eta (0))^2)^{1/2}   }\le M \qquad 
a.s.\label{res}
\ee
  This is not possible because $  \{{  \eta(t_k)-\eta(0)}/{ (E(\eta (t_{k})-  \eta (0))^2)^{1/2}   } \}$
 is a sequence of standard normal random variables.

 \medskip  Since,
   \begin{equation}
   U(0,0)-U(0,t)=E \(\( \eta (t)-  \eta (0)\) \eta (0)\)\le E \(\( \eta (t)-  \eta (0)\)^2\)^{1/2}  U^{1/2} (0,0),
   \end{equation}
  we have,
  \begin{equation}
  \frac{  U(0,0)-U(0,t)}{\phi(t)}\le \frac{  U^{1/2}(0,0) }{h(t)}.
  \end{equation} 
  Using the fact that $\lim_{t\downarrow0}h(t)=\ff$ we get  (\ref{2.30mm}).\qed

 \noindent \textbf{Proof  of Theorem \ref{theo-2.1m} }    Note that  (\ref{rev.1wqu}) implies  that 
 $\{\eta (t),t\in[0,1]\}$ and therefore $\{\eta^{2} (t),t\in[0,1]\}$ are uniformly continuous on  $[0,1]$, which in turn implies that for all $k\ge1,$ $ \{Y_{k }(t);t\in  [0,1]\}$ is uniformly continuous on  $[0,1]$.
 
  To show,
    \begin{equation}
 \lim_{h\to 0}\sup_{\stackrel{|u-v|\le h }{ u,v\in\De}}  \frac{Y_{k }(u)-Y_{k }(v) }{   \vf  (|u-v|)} \geq   2 \sup_{t\in\De}Y_{k }^{1/2}(t), \hspace{.2 in}a.s.\label{rev.2qu5}
\end{equation} 
it suffices to show that for any $d\in \De$,
\begin{equation}
 \lim_{h\to 0}\sup_{\stackrel{|u-v|\le h }{ u,v\in\De}} \frac {Y_{k }(u)-Y_{k }(v) }{    \vf  (|u-v|)} \geq  2  Y_{k }^{1/2}(d), \hspace{.2 in}a.s.\label{rev.2qu6}
\end{equation} 
This is because, (\ref{rev.2qu6}) holding almost surely implies that for any countable dense set $\De'\subset \De$,
\begin{equation}
 \lim_{h\to 0}\sup_{\stackrel{|u-v|\le h }{ u,v\in\De}}  \frac{Y_{k }(u)-Y_{k }(v) }{    \vf  (|u-v|)} \geq 2\sup_{d\in\De'}    Y_{k }^{1/2}(d), \hspace{.2 in}a.s.\label{rev.ioi}
\end{equation} 
which implies (\ref{rev.2qu5}).

Let $  u,v,d\in \De$. We write, 
\begin{eqnarray}
 \eta_{i}^{2}(u)- \eta_{i}^{2}(v)&=&( \eta_{i}(u)- \eta_{i}(v))( \eta_{i}(u)+ \eta_{i}(v))
\label{16.40u}\\
&=& ( \eta_{i}(u)- \eta_{i}(v))( 2 \eta_{i}(d)+( \eta_{i}(u)- \eta_{i}(d))+( \eta_{i}(v)- \eta_{i}(d)))   \nonumber.
\end{eqnarray}
  It follows from   (\ref{rev.1wqu}) that,											
\bea \label{2.48mm}
 && \lim_{h\to 0}\sup_{\stackrel{|u-v|\le h }{ u,v\in\De}}\frac{ \sum_{i=1}^{k}( \eta _{i }(u)-\eta _{i }(v))( \eta _{i }(u)-\eta _{i }(d)) }{     \vf (|u-v|)} \\
 &&\qquad\le  \lim_{h\to 0}\sup_{\stackrel{|u-v|\le h }{ u,v\in\De}}\frac{ \sum_{i=1}^{k}| \eta _{i }(u)-\eta _{i }(v) |\sup_{u\in\De}| \eta _{i }(u)-\eta _{i }(d)|} {    \vf (|u-v|)}\nn\\
 &&\qquad \le   \sum_{i=1}^{k}\sup_{u\in\De}| \eta _{i }(u)-\eta _{i }(d)|:=  \De^{*}.  \nn
\eea
Therefore, using (\ref{rev.2qu5})--(\ref{2.48mm}) and (\ref{2.20mm}) with 0 replaced by $d$ we see that,   
\bea
&& \lim_{h\to 0}\sup_{\stackrel{|u-v|\le h }{ u,v\in\De}}\frac{Y_{k }(u)-Y_{k }(v) }{      \vf (|u-v|)  Y^{1/2}_{k }(d)}\label{5.4mm2}\\&&\qquad \ge \lim_{h\to 0}\sup_{\stackrel{|u-v|\le h }{ u,v\in\De}}\frac{2\sum_{i=1}^{k}( \eta _{i }(u)-\eta _{i }(v)) \eta_{i}(d)}{     \|\vec\eta (d)\|_{2}  \,  \vf (|u-v|)} -\frac{2\De^{*}}{\|\vec\eta (d)\|_{2}}.\nn
   \eea
  
To simplify the notation we take $V(u,v)= {U(u,v)}/{U(d,d)}$.
Write
\begin{equation} \label{2.50mm}
  \eta_{i}(u)- \eta_{i}(v) = V(v,d) \eta_{i}(u)-V(u,d)\eta_{i}(v) +  G_{i}(u,v),
\end{equation}
where,
\begin{equation} \label{a}
  G_{i}(u,v)=(1-V(v,d) )\eta_{i}(u) -(1-V(u,d))\eta_{i}(v) .
\end{equation}
In this notation,
\begin{eqnarray}
&&\frac{\sum_{i=1}^{k}( \eta _{i }(u)-\eta _{i }(v)) \eta_{i}(d)  }{   \|\vec\eta (d)\|_{2}\,   \vf (|u-v|)}-
\frac{\sum_{i=1}^{k}G_{i}(u,v) \eta_{i}(d) }{      \|\vec\eta (d)\|_{2}\, \vf (|u-v|)}\nn\\
&&\qquad   =\frac{\sum_{i=1}^{k}\(V(v,d) \eta_{i}(u)-V(u,d)\eta_{i}(v)\) \eta_{i}(d)}{      \|\vec\eta (d)\|_{2} \,\vf (|u-v|)} .\label{2.50mm3}
\end{eqnarray}
Note that for all $u,v\in [0,1],$
\begin{equation} \nn
  E\((V(v,d) \eta_{i}(u)-V(u,d)\eta_{i}(v))\eta_{i}(d)\)=0.
\end{equation}
This shows that $\eta_{i}(d)$ is independent of $\{ V(v,d) \eta_{i}(u)-V(u,d)\eta_{i}(v);u,v\in [0,1]\}$. Therefore by  (\ref{2.32mm}),  
\bea \label{2.53mm}
 && \Big \{\sum_{i=1}^{k} (V(v,d) \eta_{i}(u)-V(u,d)\eta_{i}(v))\frac{\eta_{i}(d)}{\|\vec\eta (d)\|_{2}};\quad u,v\in [0,1]\Big\}\\
 &&\qquad \nn \stl \Big \{ V(v,d) \eta (u)-V(u,d)\eta (v);\quad u,v\in [0,1]\Big\}.
\eea
It follows that,  
\bea \label{2.54mm}
 && \lim_{h\to 0}\sup_{\stackrel{|u-v|\le h }{ u,v\in\De}}\frac{  \sum_{i=1}^{k}( \eta _{i }(u)-\eta _{i }(v) )  \eta _{i }(d) }{ \|\vec\eta (d)\|_{2}\,  \vf (|u-v|)}\\
 &&\qquad\qquad + \lim_{h\to 0}\sup_{\stackrel{|u-v|\le h }{ u,v\in\De}}\frac{\sum_{i=1}^{k}|G_{i}(u,v)| |\eta_{i}(d)| }{   \|\vec\eta (d)\|_{2} \,  \vf (|u-v|)} \nn\\&&\qquad \stackrel{law}{\geq} \nn \lim_{h\to 0}\sup_{\stackrel{|u-v|\le h }{ u,v\in\De}}\frac{   V(v,d) \eta (u)-V(u,d)\eta (v)  }{    \vf (|u-v|)}.
\eea

 Using (\ref{2.50mm})  we can write,
\bea \label{2.55mm}
 &&  \lim_{h\to 0}\sup_{\stackrel{|u-v|\le h }{ u,v\in\De}}\frac{   V(v,d) \eta (u)-V(u,d)\eta (v)   }{     \vf (|u-v|)}\\&&\qquad \ge \nn \lim_{h\to 0}\sup_{\stackrel{|u-v|\le h }{ u,v\in\De}}\frac{   (  \eta (u)- \eta (v)) }{    \vf (|u-v|)} - \lim_{h\to 0}\sup_{\stackrel{|u-v|\le h }{ u,v\in\De}}\frac{ |G(u,v)| }{     \vf (|u-v|)}\nn\\&&\qquad = 1 - \lim_{h\to 0}\sup_{\stackrel{|u-v|\le h }{ u,v\in\De}}\frac{ |G(u,v)| }{     \vf (|u-v|)}\nn,
\eea
where
\be  \label{a01}
  G (u,v)=(1-V(v,d) )\eta (u) -(1-V(u,d))\eta (v) 
  \ee 
  and we use  (\ref{rev.1wqu}) for the last line in (\ref{2.55mm}). 
  
\medskip   It follows from (\ref{2.54mm}) and (\ref{2.55mm}) that, 
   \be  \label{2.34}
    \lim_{h\to 0}\sup_{\stackrel{|u-v|\le h }{ u,v\in\De}}\frac{  \sum_{i=1}^{k}( \eta _{i }(u)-\eta _{i }(v) )  \eta _{i }(d)}{  \|\vec\eta (d)\|_{2}\,  \vf (|u-v|)}\ge 1-\HH,
\ee 
where,
\begin{equation} \nn
  \HH=  \lim_{h\to 0}\sup_{\stackrel{|u-v|\le h }{ u,v\in\De}}\frac{\sum_{i=1}^{k}|G_{i}(u,v)| |\eta_{i}(d)| }{ \|\vec\eta (d)\|_{2} \,  \vf (|u-v|)}+ \lim_{h\to 0}\sup_{\stackrel{|u-v|\le h }{ u,v\in\De}}\frac{ |G(u,v)| }{   \vf (|u-v|)}
\end{equation}

Using the Schwartz inequality followed by the triangle inequality  we note that,
\begin{equation} \nn
  \frac{\sum_{i=1}^{k}|G_{i}(u,v)| |\eta_{i}(d)|}{ \|\vec\eta (d)\|_{2} \,  \vf (|u-v|)}\le \frac{\sum_{i=1}^{k}|G_{i}(u,v)|}{\vf (|u-v|)}.
\end{equation}
Therefore, 
\begin{equation} \label{2.37mm}
  \HH\le   \sum_{i=0}^{k} \lim_{h\to 0}\sup_{\stackrel{|u-v|\le h }{ u,v\in\De}}\frac{   |G_{i}(u,v)|  }{  \vf (|u-v|)},\qquad a.s.,
\end{equation}
where for notational convenience we have set $G_{0}=G$.
 
   Set,
  \begin{equation} \label{2.34j}
 \si^{2}(u,v)=E(\eta(u)-\eta(v))^2,
\end{equation}
and  
\begin{equation}
 \wt \si^2(x)=\sup_{|u-v|\le x}\si^2(u,v).\label{ovabs.2r}
\end{equation}
 Then it follows from (\ref{rev.1wqu}) that we can write, 
 \begin{equation} \label{2.45mm}
  \vf(h)=  \wt \si(h)g(h),\mbox{ where necessarily, $\lim_{h\to 0}g(h)=\ff$}.
\end{equation}
  This follows by a minor modification of the argument used to prove (\ref{c2.12}). Note that for any sequence $h_{k}\to 0$ we can find sequences $\{u_k \}$, $\{v_k \}$ in $\De$,   with $ |u_k-v_k|\leq h_{k}$ such that $ \wt\si (h_{k})\leq 2  \si(u_k,v_k)$. Now, suppose that $\limsup_{h\to 0}g(h)=M$. Then by (\ref{rev.1wqu}) we would have, 
\begin{equation} \nn
\sup_{k\to \ff} \frac{  \eta(u_k)-\eta(v_k)}{ \si(u_k,v_k)   }\le 4M \qquad 
a.s.\label{rev.1wqumbm}
\end{equation}
This is not possible because each term $(\eta(u_k)-\eta(v_k))/ \si(u_k,v_k)  $
is a standard normal random variable.

   We show in Lemma \ref{lem-G} below that, 
  \begin{equation}
  |G(u,v)|\leq   \frac{\si( d,v )}{U^{1/2}(d,d)} |\eta (u)- \eta (v)|+\frac{\si( u,v )}{U^{1/2}(d,d)} |\eta (v)|.\label{Gboundj}
  \end{equation}
 Therefore, using (\ref{rev.1wqu}) and (\ref{2.45mm}) we have,  
\bea \label{2.60nn}
&&    \lim_{h\to 0}\sup_{\stackrel{|u-v|\le h }{ u,v\in\De}}\frac{   |G(u,v)|  }{   \vf (|u-v|)} \le    \lim_{h\to 0}\sup_{\stackrel{|u-v|\le h }{ u,v\in\De}}\frac{   |\eta (u)- \eta (v)|\si( d,v )  }{  U^{1/2}(d,d)\, \vf (|u-v|)}\,\,\,\\
  && \hspace{2in}+  \lim_{h\to 0}\sup_{\stackrel{|u-v|\le h }{ u,v\in\De}}\frac{    \si( u,v )   }{  U^{1/2}(d,d)\,  \vf (|u-v|)}\eta(v) \nn\\
  &&\qquad\qquad \le \sup_{d, v\in \De}\frac{\si( d,v )}{U^{1/2}(d,d)}+\lim_{h\to 0}\frac{1}{g(h)}\sup_{ v\in \De}\frac{|\eta\(v\)|}{U^{1/2}(d,d)} = \frac{ \wt \si(|\De|)}{U^{1/2}(d,d)}. \nn
\eea
 where $\wt \si(|\De|)$ is defined in (\ref{ovabs.2r}).
   We now use (\ref{5.4mm2}), (\ref{2.34}) and (\ref{2.60nn}) to see that,  

\bea \label{2.62mm}
      && \lim_{h\to 0}\sup_{\stackrel{|u-v|\le h }{ u,v\in\De}}\frac{ Y _{k }(u)-Y _{k }(v) }{    \vf (|u-v|)}  \\
      &&\qquad \ge 2Y^{1/2}_{k }(d)\(1- (k+1)  \frac{ \wt \si(|\De|)}{U^{1/2}(d,d)}\)-2\De^{*},\qquad a.s.,\nn
\eea
where $d\in \De$ and   $  \De^{*}$ is defined in (\ref{2.48mm}). 
  
Suppose that $|\De|=1/n$. Then it follows from 
Lemma \ref{lem-med} that,
\begin{equation} \label{2.67mm}
  P\( \De^{*}\ge k ((1+2C)  \wt \si^{2}(1/n)\log n)^{1/2}\)\le  \frac{2k}{n^{C}}.
\end{equation}
Now let $\De(d,n)\subseteq \De $, be an interval of size $1/n$ that contains $d$. It follows from (\ref{2.62mm}) and (\ref{2.67mm})
applied to      $\De(d,n) $   and  $\De^{*}(d,n) $ that  ,
\bea
  &&\label{2.68mm}  \lim_{h\to 0}\sup_{\stackrel{|u-v|\le h }{ u,v\in\De }}\frac{ Y_{k }(u)-Y_{k}(v) }{   \vf (|u-v|)}\\
  &&\qquad \geq  \lim_{h\to 0}\sup_{\stackrel{|u-v|\le h }{ u,v\in\De(d,n) }}\frac{ Y_{k }(u)-Y_{k }(v)}{   \vf (|u-v|)}\nn\\
  &&\qquad \ge 2Y^{1/2}_{k }(d)\(1-(k+1) \frac{ \wt \si(|\De|)}{U^{1/2}(d,d)}\)- 2k((1+2C)  \wt \si^{2}(1/n)\log n)^{1/2}, \nn
\eea
except, possibly, on a set of measure $2k/n^{C}$.   Taking $n\to\ff$, and using (\ref{ovabs.3}), gives (\ref{rev.2qu6}) and consequently (\ref{rev.2qu5}), which is  the lower bound in (\ref{rev.2quu}).

\medskip We now obtain  the upper bound in (\ref{rev.2quu}). Let $\De_{m,n}=\De\cap [\frac{m-1}{n}, \frac{m+1}{n}]$. 
 Analogous to (\ref{2.68mm}) we have,
\bea
  &&   \lim_{h\to 0}\sup_{\stackrel{|u-v|\le h }{ u,v\in\De_{m,n} }}\frac{ Y_{k }(u)-Y_{k }(v) }{  \vf (|u-v|)}\label{2.68pp}\\
  &&\qquad \le 2Y^{1/2}_{k }(m/n)\(1+(k+1) \frac{ \wt \si(2/n)}{U^{1/2}(\frac{m}{n},\frac{m}{n})}\)+ k((1+2C)\wt\si^{2}(2/n)\log n)^{1/2}\nn\\\nn  &&\qquad := 2 Y^{1/2}_{k }(m/n) \AA\(\wt\si(2/n),k,n,C\),\eea
except, possibly, on a set of measure $2k/n^{C}$. The proof of (\ref{2.68pp}) proceeds in essentially the same way as the proof of (\ref{2.68mm}). In the  proof of the lower bound in (\ref{2.68mm}) we subtract several terms.   In proving the upper bound in (\ref{2.68pp}) we add these terms. 
It then follows that,
\bea
 \lefteqn{   \lim_{h\to 0}\sup_{\stackrel{|u-v|\le h }{ u,v\in\De  }}\frac{ Y_{k }(u)-Y_{k }(v) }{   \vf (|u-v|)}\label{2.68pq}}\\
  &&  \le 2     \sup_{m=1,\ldots, n-1} Y^{1/2}_{k }(m/n)\AA\(\wt\si(2/n),k,n,C\)\nn\\
  && \le 2     \sup_{v\in\De } Y^{1/2}_{k }(v) \AA\(\wt\si(2/n),k,n,C\), \nn
\eea
except  possibly  on a set of measure $2k/n^{C-1}$. Taking the limit as $n\to\ff$   gives the upper bound in (\ref{rev.2quu}).
\qed
\begin{lemma}\label{lem-G}
  \begin{equation}
  |G(u,v)|\leq   \frac{\si( d,v )}{U^{1/2}(d,d)} |\eta (u)- \eta (v)|+\frac{\si( u,v )}{U^{1/2}(d,d)} |\eta (v)|.\label{Gbound}
  \end{equation}
  \end{lemma}
 {\bf  Proof: } 
 We write,
\be \label{a02}
  G (u,v)=(1-V(v,d) )(\eta (u)-\eta(v)) +(V(u,d)-V(v,d))\eta (v)\nn.
\ee
Note that,
\bea \nn
  |(V(u,d)-V(v,d))| &=&\frac{|\(E(\eta (u)- \eta (v)) \eta (d) \)|}{ U(d,d) }\\
  &\le&\frac{\(E(\eta (u)- \eta (v))^{2} E\eta^{2} (d)\)^{1/2}}{ U(d,d)}\nn\\
   &\le&\frac{\si( u,v )}{ U^{1/2}(d,d) }\nn,
\eea
and,
\bea \nn
  (1-V(v,d))  &=&\frac{\(E(\eta (d)- \eta (v)) \eta (d) \)}{ U(d,d) }\\
  &\le&\frac{\(E(\eta (d)- \eta (v))^{2} E\eta^{2} (d)\)^{1/2}}{ U(d,d)}\nn\\
   &\le&\frac{\si( d,v )}{ U^{1/2}(d,d)}  \nn.
\eea
\qed
\bl\label{lem-med}
\begin{equation} 
  P\( \sum_{i=1}^{k}\sup_{u\in\De} |\eta_i(u)-\eta_i(d)|\ge  k((1+2C) \wt \si^2(|\De|)\log 1/|\De|)^{1/2}\)\le   2k |\De|^{C}. \label{medgen}
\end{equation}
\el
{\bf  Proof: }
Let
 \textbf{a} be the median of $\sup_{u\in\De} (   \eta_1(u)-\eta_1(d))$. It follows from \cite[Lemma 5.4.3]{book} that,
\begin{equation} \label{2.64}
  P\(\sup_{u\in\De} |\eta_1(u)-\eta_1(d)|\ge \mathbf{a}+ \wt\si (|\De|)t\)\le 2e^{-t^{2}/2}.
\end{equation}
Since by \cite[(7.113)]{book},
\begin{equation} \label{2.65}
  \mathbf{a}=o(\wt\si^2(|\De|)\log 1/|\De|)^{1/2},
 \end{equation}
 we see that 
 \begin{equation} 
  P\(\sup_{u\in\De} |\eta_1(u)-\eta_1(d)|\ge  ((1+2C)\wt \si^2(|\De|)\log 1/|\De|)^{1/2}\)\le   2  |\De|^{C},\label{medgen1}
\end{equation}
which gives (\ref{medgen}).
\qed

 \begin{remark} \label{rem-2.1}{\rm When $\eta$ has stationary increments, $ \si^2(u,v)$ in (\ref{1.15}) can be written as 
 $ \si^2(u-v)$. For  $x\in [0,\de]$ define,
 \begin{equation} \nn
  \ov\si^2(x)=\mu\(s:\si^2(x)\le x\),
\end{equation}
 where $\mu$ is Lebesgue measure. Clearly,  $ \ov\si^2(x)$ is an increasing function. It is called the increasing rearrangement of $  \si^2(x)$.    (See  \cite[Section 6,4]  {book} for more details.) We show in  \cite[(6.138)]{book}
that when $\eta$ is continuous a.s.,
\begin{equation} \nn
\lim_{x\to 0}\ov\si^2(x)\log(1/x)=0.
\end{equation}
 This shows that 
  if $ \wt \si^2(x)$  is asymptotic to an increasing function at 0, then (\ref{rev.1wqu}), which implies that  $\eta$ is continuous a.s., implies (\ref{ovabs.3}).  
  
}\end{remark}

 \noindent
\begin{tabular}{lll} & \hskip20pt Michael B.  Marcus
     & \hskip20pt  Jay Rosen\\  & \hskip20pt 253 West 73rd. St., Apt. 2E
   & \hskip20pt Department of Mathematics \\    &\hskip20pt  New York, NY 10023, USA
 \hskip20pt 
     & \hskip20pt College of Staten Island, CUNY \\    & \hskip20pt mbmarcus@optonline.net
    & \hskip20pt  Staten Island, NY
10314, USA \\    & \hskip20pt      & \hskip20pt   jrosen30@optimum.net \\   & \hskip20pt 

     & \hskip20pt 
\end{tabular}

  \end{document}